\documentclass[12pt]{article}
\usepackage{amssymb,amsmath}
\usepackage{graphicx}
\usepackage{amsfonts}
\usepackage{float}
\usepackage{flafter}
\usepackage{color}

\labelwidth\leftmargini\advance\labelwidth-\labelsep

\def\R{\relax\ifmmode I\!\!R\else$I\!\!R$\fi}

\def\Z{\relax\ifmmode Z\!\!\!Z\else$Z\!\!\!Z$\fi}

\def\C{\relax\ifmmode C\!\!\!\!I\else$C\!\!\!\!I$\fi}

\def\K{\relax\ifmmode I\!\!K\else$I\!\!K$\fi}

\def\N{\relax\ifmmode I\!\!N\else$I\!\!N$\fi}

\newcounter{defcounter}[section]
{\vspace{0.1cm}\begin{sloppypar}\noindent\stepcounter{defcounter}{\bfseries
Definition
      \thesection.\thedefcounter}}%
{\end{sloppypar}\vspace{0.1cm}}

\newtheorem{theorem}{Theorem}[section]

\newtheorem{proposition}{Proposition}[section]

\newcommand{\proof}{{\bf Proof.} }

\newcommand{\qed}{\hfill $\square$}

\begin{document}
\thispagestyle{empty}
\begin{center}
{\Large {\bf Fractal attractors induced by $\beta$-shifts}}
\end{center}
\begin{center}J\"org Neunh\"auserer\\
Leibnitz University of Hannover, Germany\\
joerg.neunhaeuserer@web.de
\end{center}
\begin{center}
\begin{abstract}
We describe a class of fractal attractors induced by $\beta$-shifts. We use a coding by these shifts to show that the systems are mixing with topological entropy $\log\beta$ and have an ergodic measure of full entropy. Moreover we determine the Hausdorff dimension of the attractor.
\end{abstract}
\end{center}
\section{Introduction}
Fractal attractors are a central subject in the modern theory of dynamical systems. Famous examples coming from applied mathematics are the Lorenz attractor, the Hennon attractor, the R\"ossler attractor and the Ikeda attractor, see \cite{[HI]} for instance. Well known examples, that are of importance from a theoretical perspective, are solenoidal attractors \cite{[BO],[NE2]} and attractors of generalized Bakers maps \cite{[AY],[NE1]} .\\
We introduce here a new class of fractal attractors that are induced by $\beta$-shifts. $\beta$-shifts are intensively studied in arithmetic resp. symbolic dynamics, see \cite{[SI]} for an overview. They describe the dynamics of the expanding map $f(x)=\beta x\mbox{ mod }1$ on the Intervall $[0,1)$, where $\beta>1$. ~\\
For parameters $\beta\in(1,2)$ and $\tau\in(0,0.5)$ let us consider the map $f:[0,1]^2\to [0,1]^2$ given by
\[
     f(x,y)=\left\{\begin{array}{ll} (\beta x,\tau y), & x\in [0,\beta^{-1}] \\
         (\beta x-1,\tau y+(1-\tau)), & x\in(\beta^{-1},1]\end{array}\right..  \]
We define the compact attractor of the dynamical system $([0,1]^2,f)$ by
\[ \Lambda_{\beta,\tau}=\mbox{closure}(\bigcap_{i=0}^{\infty}f^i([0,1]^2)). \]
\begin{center}
\begin{figure}
\vspace{0pt}\hspace{0pt}\scalebox{0.5}{\includegraphics
{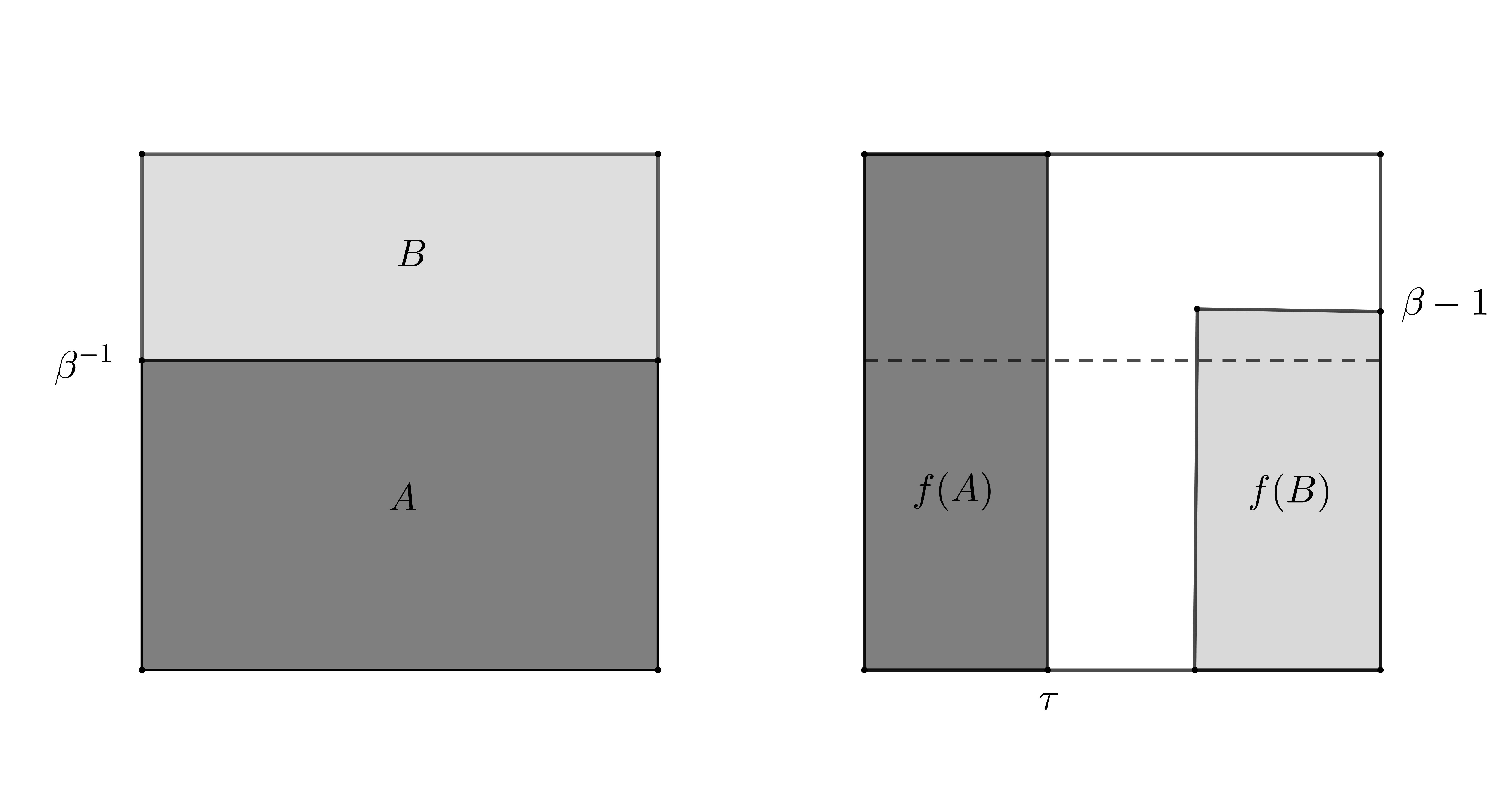}}
\caption{The action of $f$ on $[0,1]^2$}
\end{figure}
\end{center}
In the overlapping case $\tau\in(0.5,1)$ this attractor was studied in \cite{[FP]}, especially the question if there is an absolutely continuous ergodic measure for the system is addressed. We consider here the non-overlapping case $\tau\in(0,0.5)$. ~\\
In the next section we will describe the dynamics of $f$ on the attractor $\Lambda_{\beta,\tau}$ symbolically using $\beta$-shifts. With the help of this  description we will show that the system $(\Lambda_{\beta,\tau},f)$ is topological mixing.
In section three we study the entropy of the dynamical system. Using the symbolic coding we show that topological entropy of the system is $\log(\beta)$ and that there is an ergodic measure of full entropy. In the last section we determine the Hausdorff dimension of the attractor $\Lambda_{\beta,\tau}$ which turn to be in $(1,2)$ for all $\beta\in(1,2)$ and $\tau\in(0,0.5)$. This means that the attractor is in fact a fractal according to the usual definition.
\section{Symbolic dynamics}
Consider the space of bi-infinite sequences on two symbols $\Sigma=\{0,1\}^{\mathbb{Z}}$ with the natural product topology which is induced by the metric
\[ d((s_{k}),(t_{k}))=\sum_{i=0}^{\infty}|s_{k}-t_{k}|2^{-|k|}.\]
The shift $\sigma:\Sigma\to\Sigma$ given by $\sigma((s_{k}))=(s_{k-1})$ is an universal model in chaotic dynamics, see \cite{[KH]} for instance.
For a real number $\beta\in(1,2)$ we consider here a subshift given by
\[ X_{\beta}=\{(s_{k})\in\Sigma~|~\sum_{k=1}^{\infty}s_{i-k}\beta^{-k}<1\mbox{ for all }i\in\mathbb{Z}\}\]
and its closure
\[\overline{X}_{\beta}=\{(s_{k})\in\Sigma~|~\sum_{k=1}^{\infty}s_{i-k}\beta^{-k}\le 1\mbox{ for all }i\in\mathbb{Z}\}.\]
In addition we use
\[\overline{X}_{\beta}^{\star}=\overline{X}_{\beta}\backslash\{ (s_{k})|\exists i\in\mathbb{Z}:s_{i-k}=0\mbox{ for all }k\in\mathbb{Z}\}.\]
The sets $X_{\beta}$, $\overline{X}_{\beta}$ and $\overline{X}_{\beta}^{\star}$ are obviously forward and backward invariant under the shift maps $\sigma$. The dynamical systems $(\overline{X}_{\beta},\sigma)$ are know as two-sided $\beta$-shift.~\\
 Now we introduce a coding map $\pi:\overline{X}_{\beta}\to Q$ by
\[\pi((s_{k}))=(\sum_{k=1}^{\infty}s_{-k}\beta^{-k},\sum_{k=0}^{\infty}s_{k}(1-\tau)\tau^{k}).\]
This map has the following properties:
\begin{proposition}
$\pi$ is continuous with $\pi(\overline{X_{\beta}})=\Lambda_{\beta,\tau}$ and the map conjugates $f$ and the shift $\sigma$ on $\overline{X}_{\beta}^{\star}$, that means
\[ f(\pi((s_{k})))=\pi(\sigma((s_{k})))\]
for all sequences in $(s_{k})\in \overline{X_{\beta}}^{\star}$. Moreover $f$ is injective on $X_{\beta}$.
\end{proposition}
\proof Let $(s_{k}^{(i)})$ be a sequence of sequences in $\overline{X_{\beta}}$ with $\lim_{i\to\infty}(s_{k}^{(i)})=(s_{k})$. Let
\[ M(i)=\max\{n~|~s_{k}^{(i)}=s_{k},~k\in\{-n,\dots,-1,0,1,\dots, n\}\}\]
By the definition of the metric on $\overline{X_{\beta}}$ we get $\lim_{i\to\infty}M(i)=\infty$. Looking at the definition of the coding map $\pi$ this obviously implies $\lim_{i\to\infty}\pi((s_{k}^{(i)}))=\pi((s_{k}))$. Hence $\pi$ is continuous. We now show the conjugacy. Consider a sequence $(s_{k})\in \overline{X_{\beta}}^{\star}$. If $\sum_{k=1}^{\infty}s_{-k}\beta^{-k}\le\beta^{-1}$ we have $s_{-1}=0$ ($s_{-1}=1$ would imply $s_{-k}=0$ for all $k\ge2$). Hence we get
\[ f(\pi((s_{k})))=(\beta\sum_{k=1}^{\infty}s_{-k}\beta^{-k},\tau\sum_{k=0}^{\infty}s_{k}(1-\tau)\tau^{k})\]
\[ =(\sum_{k=1}^{\infty}s_{-k-1}\beta^{-k},\sum_{k=0}^{\infty}s_{k-1}(1-\tau)\tau^{k})=\pi(\sigma((s_{k}))). \]
If $\sum_{k=1}^{\infty}s_{-k}\beta^{-k}> \beta^{-1}$ we have $s_{-1}=1$, an thus
\[ f(\pi((s_{k})))=(\beta\sum_{k=1}^{\infty}s_{-k}\beta^{-k}-1,\tau\sum_{k=0}^{\infty}s_{k}(1-\tau)\tau^{k}+(1-\tau))\]\[=
(\sum_{k=1}^{\infty}s_{-k-1}\beta^{-k},\sum_{k=0}^{\infty}s_{k-1}(1-\tau)\tau^{k})=\pi(\sigma((s_{k}))).\]
Now consider the map $f$ on $X_{\beta}$ and let $(s_{k}),(t_{k})\in X_{\beta}$. By the definition of $X_{\beta}$
\[ \sum_{k=1}^{\infty}s_{-k}\beta^{-k}=\sum_{k=1}^{\infty}t_{-k}\beta^{-k}\]
implies $s_{-k}=t_{-k}$ for all $k\in\mathbb{N}$. Since $\tau<1/2$
\[ \sum_{k=0}^{\infty}s_{k}(1-\tau)\tau^{k}=\sum_{k=0}^{\infty}t_{k}(1-\tau)\tau^{k}\]
implies  $s_{k}=t_{k}$ for all $k\in\mathbb{N}_{0}$. Hence $f$ is injective on $X_{\beta}$.
We now define
\[ S_{\beta}=\{(s_{k})\in\Sigma~|~\sum_{k=1}^{\infty}s_{i-k}\beta^{-k}\le 1\mbox{ for all }i\le 0\}\backslash\{ (s_{k})|\exists i\in\mathbb{N}:s_{i-k}=0\mbox{ for all }k\in\mathbb{Z}\}. \]
Note that $\bigcap_{i=0}^{\infty}\sigma^i( S_{\beta})=\overline{X}_{\beta}^{\star}$.
Let
\[ C_{\tau}= \{\sum_{k=0}^{\infty}s_{k}(1-\tau)\tau^{k}~|~s_{k}\in\{0,1\},k\in\mathbb{N}_{0}\}\]
and
\[ I_{\beta}=[0,1]\backslash \{\sum_{k=1}^{n}s_{k}\beta^{-k}~|~s_{k}\in\{0,1\},k=1,\dots,n\}\]
We have
\[        I_{\beta}\times C \subseteq \pi(S_{\beta})\subseteq [0,1]\times C.\]
Since $f(\pi(S_{\beta}))=\pi(\sigma(S_{\beta})$ we get
\[ f^{n}(I_{\beta}\times C)\subseteq \pi(\sigma^n(S_{\beta}))\subseteq f^{n}([0,1]\times C)\]
and hence
\[  \bigcap_{n=0}^{\infty}f^{n}(I_{\beta}\times C)\subseteq \pi(\overline{X}_{\beta}^{\star})\subseteq \bigcap_{n=0}^{\infty}f^{n}([0,1]\times C).\]
The closure of $I_{\beta}$ is $[0,1]$ and the closure of $\overline{X}_{\beta}^{\star}$ is $\overline{X}_{\beta}$. Thus we obtain $\pi(\overline{X_{\beta}})=\Lambda_{\beta,\tau}$.
  \qed~\\~\\
Using this proposition we get the following result on the dynamics of $f$ on the attractor $\Lambda_{\beta,\tau}$.
\begin{theorem}
The dynamical system $(\Lambda_{\beta,\tau},f)$ is topological mixing:  For nonempty open sets $A$ and $B$ there exists an integer $N$ such that, for all $n > N$
\[ f^{n}(A)\cap B\not=\emptyset.\]
\end{theorem}
\proof It is known that the system $(\overline{X_{\beta}},\sigma)$ is topological mixing, see \cite{[RE]}. Let $A$ and $B$ be two open sets in $\Lambda_{\beta,\tau}$. By continuity of $\pi$ the preimages  $\pi^{-1}(A)$ and $\pi^{-1}(B)$ are open and they are contained in $\overline{X_{\beta}}$. Hence there ist an $N$, such that for all $n>N$ there exists a sequence $s(n)\in \overline{X}_{\beta}$, such that
\[s(n)\in \sigma^{n}(\pi^{-1}(A))\cap\pi^{-1}(B). \]
Since the set $\sigma^{n}(\pi^{-1}(A))\cap\pi^{-1}(B)$ is open in $\overline{X}_{\beta}$ we may assume that $s(n)\in \overline{X}_{\beta}^{\star}$.
Hence we get:
 \[\pi(s(n))\in \pi(\sigma^{n}(\pi^{-1}(A))\cap\pi^{-1}(B))\subseteq \pi(\sigma^{n}(\pi^{-1}(A)))\cap\pi(\pi^{-1}(B)) \]
\[ =f^{n}(\pi((\pi^{-1}(A)))\cap\pi(\pi^{-1}(B))\subseteq f^{n}(A)\cap B\not=\emptyset.\]
Thus $(\Lambda_{\beta,\tau},f)$ is topological mixing.\qed
\section{Entropy}
For convenience we first recall the definition of the topological and metric entropy. Let $(X,f)$ be dynamical system on a compact space $X$. The entropy of an open covering ${\mathfrak U}$ of $X$ is given $H({\mathfrak U})=\log\sharp{\mathfrak U}$, where $\sharp{\mathfrak U}$ is the minimal number of elements in ${\mathfrak U}$  that cover $X$.
The entropy  of the system $(X,f)$ with respect to ${\mathfrak U}$ is
\[ h(f,{\mathfrak U})=\lim_{n\to\infty}\frac{1}{n}H({\mathfrak U}\vee f^{-1}({\mathfrak U})\vee\dots\vee f^{-n}({\mathfrak U})),\]
where a covering ${\mathfrak U_{1}}\vee {\mathfrak U_{2}}$ consists of the intersections of elements in ${\mathfrak U}_{1}$ and ${\mathfrak U}_{2}$.
The topological entropy of the system is
\[ h(f)=\sup\{h(f,{\mathfrak U})~|~{\mathfrak U}\mbox{ is an open covering of }X\}.\]
A Borel probability measure $\mu$ on $X$ is ergodic with respect to $f$ if it is invariant, $\mu\circ f^{-1}=\mu$, and sets $B$ with $f^{-1}(B)=B$ have measure zero or one. The metric entropy of a system $(X,f,\mu)$ with respect to a measurable partition ${\mathfrak P}$ is
\[ h(f,\mu,{\mathfrak P})=\lim_{n\to\infty}\frac{1}{n}H(\mu,{\mathfrak P}\vee f^{-1}({\mathfrak P})\vee\dots\vee f^{-n}({\mathfrak P})),\]
where is entropy of a measurable partition is
\[ H(\mu,\mathfrak{ P})=-\sum_{P\in\mathfrak{P}}\mu(P)\log(\mu(P)).\]
The metric entropy of the system is
\[ h(f,\mu)=\sup\{h(f,{\mathfrak P})~|~{\mathfrak P}\mbox{ is a measurable partition of }X\}.\]
We recommend \cite{[WA]} and \cite{[KH]} for an introduction to entropy theory. ~\\~\\
The entropy of $\beta$-shifts is well studied. The topological entropy of $(\overline{X_{\beta}},\sigma)$ is given by $\log\beta$, $h(\sigma_{| \overline{X_{\beta}}})=\log\beta$, see \cite{[RE]}. Moreover there is a shift ergodic measure $\mu_{\beta}$ on $\overline{X_{\beta}}$ with full entropy $h(\mu_{\beta},\sigma)=\log\beta$, see \cite{[PA]}. We will transfer these results to the dynamical system $(\Lambda_{\beta,\tau},f)$ using the symbolic coding in proposition 2.1. and prove:
\begin{theorem}
The dynamical system $(\Lambda_{\beta,\tau},f)$ has topological entropy $\log(\beta)$ and there is an ergodic measure of full entropy.
\end{theorem}
\proof Proposition 2.1. shows in the terminology of topological dynamics that $(\Lambda_{\beta,\tau},f)$ is a factor of $(\overline{X_{\beta}}^{\star},\sigma)$. It is well known and straightforward from the definition above, that this implies
\[ h(f_{|\Lambda_{\beta,\tau}})\le h(\sigma_{| \overline{X_{\beta}}^{\star}}) \le h(\sigma_{| \overline{X_{\beta}}})=\log\beta.\]
Let $\mu$ be an ergodic measure for $(\overline{X_{\beta}},\sigma)$. We show first that $\mu(\overline{X_{\beta}}\backslash \overline{X_{\beta}}^{\star})=0$. Let \[ A_{k}=\{(s_{i})~|~s_{k}=1,~s_{i}=0\mbox{ for }i<k\}\] for $k\in\mathbb{Z}$. The shifted sets $\sigma^i(A_{k})$ and $\sigma^j(A_{k})$ are disjoint for $i,j\in \mathbb{N}_{0}$ with $i\not=j$. Since $\mu$ ist $\sigma$-invariant this implies $\mu(A_{k})=0$, so $\mu(\{ (s_{k})|\exists i\in\mathbb{Z}:s_{i-k}=0\mbox{ for all }k\in\mathbb{Z}\})=0$.~\\
Now we prove $\mu(\overline{X_{\beta}}\backslash X_{\beta})=0$.
We may decompose $\overline{X_{\beta}}\backslash X_{\beta}$ in the following way
\[ \overline{X_{\beta}}\backslash X_{\beta}=\{(s_{k})\in\Sigma~|\mbox{For some }i\in\mathbb{Z}:~\sum_{k=1}^{\infty}s_{i-k}\beta^{-k}=1\mbox{ and }s_{k}=0,~k\ge i\}=\bigcup_{i=-\infty}^{\infty} N_{i},\]
where $N_{i}$ contains all sequences $(s_{k})\in \overline{X_{\beta}}\backslash X_{\beta}$ with $s_{i-1}=1$ and $s_{k}=0$ for $k\ge i$. Note that $\sigma^{-a}(N_{i})\cap \sigma^{-b}(N_{i})=\emptyset$ for $a\not=b$. Since $\mu(\sigma^{-a}(N_{i}))=\mu(N_{i})$ this implies $\mu(N_{i})=0$ and $\mu(\overline{X_{\beta}}\backslash X_{\beta})=0$.\\
Now we project $\mu$ to $\Lambda_{\beta,\tau}$ via $\nu=\pi(\mu)=\mu\circ \pi^{-1}$. By proportion 2.1 and the consideration above the measure space $(\overline{X_{\beta}},\mu)$ and $(\Lambda_{\beta,\tau},\nu)$ are measure theoretical isomorphic. Moreover the dynamical systems
$(\overline{X_{\beta}},f,\mu)$ and $(\Lambda_{\beta,\tau},\sigma,\nu)$ are by proposition 2.1 measure theoretical conjugated with $f\circ \pi=\pi \circ \sigma$. It is well known and straightforward from the definition above, that this implies \[ h(f_{|\Lambda_{\beta,\tau}},\nu)=h(\sigma_{| \overline{X_{\beta}}},\mu).\]
If $\mu_{\beta}$ is the measure of full entropy for the $\beta$-shift the projected measure $\nu_{\beta}$ has full entropy for $f$,  \[ h(f_{|\Lambda_{\beta,\tau}},\nu_{\beta})=\log\beta= h(f_{|\Lambda_{\beta,\tau}}).\]
Here we use the fact that the metric entropy is always bounded by topological entropy of a dynamical system.
\qed
\section{Dimension} In this section we determine the Hausdorff dimension of the attractor $\Lambda_{\beta,\tau}$ defined in section one. We refer to \cite{[FA]} or \cite{[PE]} for an introduction to dimension theory.  Let us recall that the $d$-dimensional Hausdorff measure of a set $B\subseteq \mathbb{R}^{2}$ is given by
 \[ H^{d}(B)=\lim_{\epsilon \longmapsto 0}\inf\{\sum_{i=1}^{\infty} \mbox{diameter}(C_{i})^{d}|B\subseteq \bigcup_{i=1}^{\infty}C_{i},~\mbox{diameter}(C_{i})<\epsilon\}\]
and the Hausdorff dimension of $B$ is
\[ \dim_{H}B=\inf\{d~|~H^{d}(B)=0\}=\sup\{d~|~H^{d}(B)=\infty\}.\]
As an upper bound on the Hausdorff dimension we will use the (lower) box-counting dimension:
\[ \dim_{H}B\le \underline{\lim}_{\epsilon\to 0}\frac{\log N_{\epsilon}(B)}{\log\epsilon^{-1}},\]
 where $N_{\epsilon}(B)$ is the minimal number of squares with side length $\epsilon$ needed to cover $B$. As a lower bound we will use the Hausdorff dimension of a Borel probability measure $\mu$, which is given by
 \[  \dim_{H}\mu=\inf\{\dim_{H}B~|~\mu(B)=1\}. \]
 We prove the following result:
 \begin{theorem} For all $\beta\in(1,2)$ and $\tau\in(0,0.5)$ we have
\[ \dim_{H}\Lambda_{\beta,\tau}=1+\frac{\log\beta}{\log\tau^{-1}}\]
 \end{theorem}
 \proof If $R$ is an aligned rectangle in $[-1,1]^2$, than $f(R)$ consists of one or two aligned rectangles. For simplicity we refer to a line segment here as a rectangle of side length zero. $f^{n}([-1,1]^2)$ consists of at most $2^n$ aligned rectangles $R_{1},R_{2},\dots,R_{t}$. In the first coordinate direction $f$ is an expansion with factor $\beta$ hence $\sum_{i=1}^{t}x_{i}=\beta^{n}$, where $x_{i}$ is the length of $R_{i}$ in the first coordinate direction. In the second coordinate direction $f$ is a contraction with factor $\tau$. The length of $R_{i}$ in the second coordinate direction is $\tau^{n}$. The number of squares of side length $\tau^{n}$, needed to cover $R_{i}$, is less than $(x_{i}/\tau^{n}+1)$. Hence we have
 \[ N_{\tau^{n}}(\Lambda_{\beta,\tau})\le  N_{\tau^{n}}(f^{n}([0,1]^2))\le \beta^{n}/\tau^{n}+t \le \beta^{n}/\tau^{n}+2^n\le(\beta^{n}+1)/\tau^{n}\]
and obtain
\[ \dim_{H}\Lambda_{\beta,\tau}\le \lim_{n\to\infty}\frac{\log N_{\tau^{n}}(\Lambda_{\beta,\tau})}{\log(\tau^{-n})}=\lim_{n\to\infty}\frac{\log((\beta^{n}+1)/\tau^{n})}{\log(\tau^{-n})}= 1+\frac{\log\beta}{\log\tau^{-1}}.\]
For a $f$-ergodic measure $\nu$ on $\Lambda_{\beta,\tau}$ we have the Ledrappier-Young formula for the dimension of the measure
\[ \dim_{H}\nu=h(f,\nu)(\frac{1}{\log\beta}+\frac{1}{\log\tau^-1}),\]
see \cite{[LY],[Y]}. The theory of Ledappier-Young is formulated for differentiable systems without singularity, but it may be applied in our context as well. The argument for this fact is given in \cite{[NE2]}, one has to guarantee existence of Lyapunov charts. Let $\nu_{\beta}$ now be the ergodic measure of full entropy for $f$ described in the last section. We have
\[ \dim_{H}\nu_{\beta}=1+\frac{\log\beta}{\log\tau^{-1}}\ge\dim_{H}\Lambda_{\beta,\tau}, \]
which completes the proof.\qed

\end{document}